\documentclass[reqno,11pt]{amsart}
\usepackage{amsfonts,amsmath,amssymb,amsthm,amscd}
\usepackage{mathrsfs}
\usepackage{upgreek}
\usepackage{color}
\usepackage{graphicx}
\usepackage{calc}
\usepackage{setspace}
\usepackage{amscd}
\usepackage[mathscr]{eucal}
\DeclareFontFamily{OT1}{pzc}{}
\DeclareFontShape{OT1}{pzc}{m}{it}{<-> s * [1.050] pzcmi7t}{}
\DeclareMathAlphabet{\mathpzc}{OT1}{pzc}{m}{it}
\newenvironment{frcseries}{\fontfamily{frc}\selectfont}{}

\newlength{\wcwidth}
\newlength{\wcheight}
\newcommand{\widecheck}[1]{\ensuremath{
\settowidth{\wcwidth}{#1}
\settoheight{\wcheight}{#1}
\addtolength{\wcheight}{1pt}
\makebox[0.0cm][l]{%
\raisebox{\depth+\wcheight}[0cm][0cm]{%
\hspace{-0.025in}\scalebox{-1}{$\widehat{\hphantom{#1}}$}}}#1
\rule{0pt}{\wcheight +2.5pt}}}
\newlength{\wcswidth}
\newlength{\wcsheight}

\usepackage{hyperref}
\hypersetup{backref,colorlinks=true,citecolor=blue,linkcolor=blue}
\title[Heegaard Floer homology and Seiberg--Witten Floer homology]{HF$=$HM I : Heegaard Floer homology\\ and\\ Seiberg--Witten Floer homology}
\author[Kutluhan]{\c{C}a\u{g}atay Kutluhan}
\author[Lee]{Yi-Jen Lee}
\author[Taubes]{Clifford Henry Taubes}
\address{Department of Mathematics, Columbia University, New York, NY 10027}
\email{kutluhan@math.columbia.edu}
\address{Department of Mathematics, Purdue University, West Lafayette, IN 47907}
\email{yjlee@math.purdue.edu}
\address{Department of Mathematics, Harvard University, Cambridge, MA 02138}
\email{chtaubes@math.harvard.edu}
\newtheorem{theorem}{Theorem}[section]

\newtheorem{lemma}[theorem]{Lemma}

\newtheorem*{mt}{Main Theorem}

\theoremstyle{definition}

\newtheorem{remark}{Remark}
\newtheorem*{rmk}{Remark}
\newtheorem{property}{Property}

\newtheorem*{ack}{Acknowledgments}

\numberwithin{equation}{section}

\newcommand{\m}{\mathrm{M}}
\newcommand{\y}{\mathrm{Y}}
\newcommand{\ym}{\overline{\mathrm{Y}}}
\newcommand{\yp}{\mathrm{Y}}
\newcommand{\umap}{\mathbb{U}}

\newcommand{\sphere}{\mathrm{S}^2}
\newcommand{\handle}{\mathrm{S^1\times S^2}}
\newcommand{\canon}{{\mathrm{K}^{-1}}}

\newcommand{\pair}{\mathfrak{p}}
\newcommand{\ho}{\mathcal{H}_{0}}
\newcommand{\hp}{\mathcal{H}_{\pair}}
\newcommand{\z}{\mathcal{Z}}

\newcommand{\zhf}{\z_{\mathrm{HF}}}
\newcommand{\delhf}{\partial_{\mathrm{HF}}}
\newcommand{\zsw}{\z_{\mathrm{SW},\y,r}}
\newcommand{\zhsw}{\hat{\z}_{\mathrm{SW},\y,r}}
\newcommand{\delsw}{\partial_{\mathrm{SW}}}
\newcommand{\diff}{\mathrm{d}}

\newcommand{\gu}{\mathrm{u}}
\newcommand{\psg}{\mathfrak{v}}
\newcommand{\go}{\mathcal{G}_{\m_{\Lambda}}}
\newcommand{\hfplus}{\mathrm{HF}^{\hspace{0.0125in}+\hspace{0.0125in}}}
\newcommand{\hfminus}{\mathrm{HF}^{\hspace{0.0125in}-\hspace{0.0125in}}}
\newcommand{\hfinfty}{\mathrm{HF}^\infty}
\newcommand{\hfhat}{\widehat{\mathrm{HF}}}

\newcommand{\hmcheck}{\mathrm{\widecheck{HM}}_\ast}
\newcommand{\hmhat}{\widehat{\mathrm{HM}}_\ast}
\newcommand{\hmbar}{\overline{\mathrm{HM}}_\ast}
\newcommand{\hm}{\widetilde{\mathrm{HM}}}

\newcommand{\hhcheck}{\mathrm{H}_\ast^+}
\newcommand{\hhhat}{\mathrm{H}_\ast^-}
\newcommand{\hhbar}{\mathrm{H}_\ast^\infty}

\newcommand{\ech}{\mathit{ech}}
\newcommand{\echplus}{\mathit{ech}^{\hspace{0.0125in}+\hspace{0.0125in}}}
\newcommand{\echminus}{\mathit{ech}^{\hspace{0.0125in}-\hspace{0.0125in}}}
\newcommand{\echinfty}{\mathit{ech}^\infty}

\newcommand{\vhat}{\hat{\mathrm{V}}}
\newcommand{\spinc}{\mathrm{Spin^{\mathbb{C}}}}
\begin{document}
\begin{abstract}
Let $\mathrm{M}$ be a closed, connected, and oriented $3$-manifold. This article is the first of a five part series that constructs an isomorphism between the Heegaard Floer homology groups of $\m$ and the corresponding Seiberg--Witten Floer homology groups of $\m$. \end{abstract}
\maketitle
\section{Introduction}
\label{s1}
\par This article and its sequels describe an isomorphism between the Heegaard Floer homology of a given closed, connected and oriented 3-manifold and the balanced version of its Seiberg--Witten Floer homology.  This article gives an overview of the proof, leaving all but a few of the technical details to the sequels.  What follows directly sets the stage for a formal statement of the equivalence.
\par Let $\m$ be a closed, connected, and oriented $3$-manifold. As explained in the book \cite{km} by P. B. Kronheimer and T. S. Mrowka, given a $\spinc$ structure $\mathfrak{s}$ on $\m$ there are three different flavors of balanced Seiberg--Witten Floer homology groups with coefficient ring $\mathbb{Z}$. These are denoted by $\hmbar(\m,\mathfrak{s},c_b)$, $\hmhat(\m,\mathfrak{s},c_b)$ and $\hmcheck(\m,\mathfrak{s},c_b)$. Each of these groups is endowed with a relative grading by a certain quotient of the group $\mathbb{Z}$ determined by the given $\spinc$ structure, and a canonical $\mathbb{Z}[\umap]\otimes\wedge^\ast(\mathrm{H}_1(\m;\mathbb{Z})/tor)$-module structure. Moreover, these groups fit into a long exact sequence where the homomorphisms respect the $\mathbb{Z}[\umap]\otimes\wedge^\ast(\mathrm{H}_1(\m;\mathbb{Z})/tor)$-module structures. 
\par The Heegaard Floer homology groups of $\m$ are graded Abelian groups defined by P. Ozsv\'ath and Z. Szab\'o in \cite{os1} and \cite{os2}.  These groups are also labeled by $\spinc$ structures on $\m$, and given a $\spinc$ structure $\mathfrak{s}$ on $\m$ there are three different flavors of Heegaard Floer homology groups. These are denoted by $\hfinfty(\m,\mathfrak{s})$, $\hfminus(\m,\mathfrak{s})$, and $\hfplus(\m,\mathfrak{s})$.  Each of these groups admits a relative grading by the same group as its Seiberg--Witten counterpart and a canonical $\mathbb{Z}[\umap]\otimes\wedge^\ast(\mathrm{H}_1(\m;\mathbb{Z})/tor)$-module structure. Furthermore, these groups also fit into a long exact sequence where the homomorphisms respect the $\mathbb{Z}[\umap]\otimes\wedge^\ast(\mathrm{H}_1(\m;\mathbb{Z})/tor)$-module structures. 
\par With the preceding understood, this article and its sequels \cite{klt1}--\cite{klt4} prove the following theorem.
\begin{mt}
Fix a $Spin^{\mathbb{C}}$ structure $\mathfrak{s}$ on $\m$. There exists a commutative diagram
\begin{equation*}
\begin{CD}
\llap{\(\cdots  \)}\hfminus(\m,\mathfrak{s}) @>>> \hfinfty(\m,\mathfrak{s})
@>>>\hfplus(\m,\mathfrak{s}) \rlap{\(\cdots \)}\\
@VVV @VVV @VVV \\
\llap{\(\cdots \)}\hmhat(\m,\mathfrak{s},c_b) @>>>
\hmbar(\m,\mathfrak{s},c_b)@>>>
\hmcheck(\m,\mathfrak{s},c_b)\rlap{\(\cdots  \)}
\end{CD}
\end{equation*}
where the vertical arrows are isomorphisms and the top and bottom rows are the respective long exact sequences for Heegaard Floer homology and Seiberg--Witten Floer homology. The vertical homomorphisms preserve the relative gradings and intertwine the respective $\mathbb{Z}[\umap]\otimes\wedge^\ast(\mathrm{H}_1(\m;\mathbb{Z})/tor)$-module structures.
\end{mt}

Our proof of the Main Theorem invokes a third sort of homology theory for $3$-manifolds. This is a version of M. Hutchings's \emph{embedded contact homology} (see \cite{hu3}). The embedded contact homology groups in the present context are defined with the choice of a \emph{stable Hamiltonian structure}. A stable Hamiltonian structure on a closed, oriented $3$-manifold is a pair $(a,w)$ where $a$ is a $1$-form and $w$ is a closed $2$-form and $\diff a $ is in the span of $w$. The embedded contact homology groups are also labeled in part by the $\spinc$ structures on the ambient manifold. Moreover, they admit relative gradings that are analogous to those of the Heegaard Floer homology groups, and an analogous module structure for the tensor product of $\mathbb{Z}[\umap]$ with the exterior product of the first homology modulo torsion. The version used here is defined for a particular stable Hamiltonian structure and for certain $\spinc$ structures on the connected sum of $\m$ with certain number of copies of $\handle$. In what follows, we denote by $\yp$ the connected sum with the orientation induced from $\m$ and by $\ym$ the connected sum with the orientation induced from $-\m$. Of interest here is the manifold $\ym$. The details of our construction of $\ym$ and the particular stable Hamiltonian structure are given in Section 1 of \cite{klt1}. It suffices to say here that $\ym$ and its geometry are constructed using a chosen $\spinc$ structure on $\m$, and the data that is used to define the corresponding Heegaard Floer homology groups of $\m$. Use $\mathfrak{s}$ to denote this $\spinc$ structure. 
\par We consider a variant of embedded contact homology on $\ym$ that is defined for a $\spinc$ structure that corresponds in a natural way to the chosen $\spinc$ structure $\mathfrak{s}$. This variant is denoted by $\echinfty$ and it is of the sort described by Definition 11.8 in \cite{hs}. The special geometry of $\ym$ is used to write the $\echinfty$ chain complex and the differential using the chain complex and the differential for $\hfinfty(\m,\mathfrak{s})$. We say more about this in Section~\ref{s2}. Theorems \ref{t2.2}--\ref{t2.4} describe the relationship between the relevant versions of embedded contact homology of $\ym$ and the Heegaard Floer homology of $\m$.
\par We also use the special geometry of $\ym$ to identify our $\echinfty$ chain complex with a chain complex that computes the three flavors of a version of the Seiberg--Witten Floer homology on $\yp$ with a particular local coefficient system. This is a version of the sort of Seiberg--Witten Floer homology that is described in Section 3.7 of \cite{km}. The relevant version of the Seiberg--Witten equations on $\yp$ and the chain complex that computes the corresponding three homology groups are described in Section~\ref{s3}. Theorem \ref{t3.3} states the relationship between the Seiberg--Witten Floer chain complex and the chain complex for $\echinfty$. Theorem \ref{t3.4} exploits Theorem \ref{t3.3} to identify the relevant Seiberg--Witten Floer homology groups on $\yp$ with three homology groups computed using the embedded contact homology chain complex on $\ym$. 
\par The identifications given by Theorems \ref{t2.4} and \ref{t3.4} between our version of embedded contact homology on $\ym$ and the Heegaard Floer homology on $\m$ and our version of Seiberg--Witten Floer homology on $\yp$ are used to write the latter groups in terms of the former. Meanwhile, a connected sum formula is used to write these same Seiberg--Witten Floer homology groups in terms of the groups $\hmbar(\m,\mathfrak{s},c_b)$, $\hmhat(\m,\mathfrak{s},c_b)$ and $\hmcheck(\m,\mathfrak{s},c_b)$. The form of the connected sum formula is analogous to a corresponding Heegaard Floer connected sum formula in \cite{os1} and was suggested to the second author by Mrowka and Ozsv\'{a}th. The connected sum formula is stated in Section \ref{s4}. These identifications lead directly to a proof of the Main Theorem. In particular, our Main Theorem follows directly from the upcoming Theorems \ref{t2.4}, \ref{t3.4}, and \ref{t4.1}.
\par By way of summary, our proof of the Main Theorem involves establishing the following correspondences:
\vspace{11pt}
\begin{center}
\begin{math}
\begin{array}{ccc}
\mathrm{SW\;Floer\;homology\;of\;\yp} & \leftrightarrow & \mathrm{embedded\;contact\;homology\;of\;\ym}\\
\updownarrow & & \updownarrow \\
\mathrm{SW\;Floer\;homology\;of\;\m} & \simeq & \mathrm{Heegaard\;Floer\;homology\;of\;\m}. 
\end{array}
\end{math}
\end{center}
\vspace{11pt}
\par There is also a fourth flavor of the Seiberg--Witten Floer homology of $\m$ defined by the second author in \cite{yjl}. This is denoted here by $\hm(\m,\mathfrak{s})$. In the meantime, Ozsv\'{a}th and Szab\'{o} define in \cite{os1} a fourth version of Heegaard Floer homology, denoted by $\hfhat(\m,\mathfrak{s})$. It follows from our proof of the Main Theorem that $\hm(\m,\mathfrak{s})$ is isomorphic to $\hfhat(\m,\mathfrak{s})$. V. Colin, P. Ghiggini and K. Honda have recently announced another proof of this isomorphism using open book decompositions and Hutchings's embedded contact homology for adapted contact forms.
\par The notion of using connected sums of $\m$ with $\handle$ to relate the Heegaard Floer homology on $\m$ with some version of Seiberg--Witten Floer homology has antecedents in the work of the second author \cite{yjl} on an as yet unsuccessful strategy to prove the Main Theorem. The appearance of the local coefficient version of Seiberg--Witten Floer homology was foreseen and a version of the connected sum formula stated in Section \ref{s4} was proved in a somewhat different context by the second author in \cite{yjl}. Bloom, Mrowka and Ozsv\'{a}th \cite{bmo} have proved a related connected sum formula for applications in a more general context.
\par The relationship between $\echinfty$ on $\ym$ and a version of the Seiberg--Witten Floer homology on $\yp$ is the analog here of the theorem by the third author proved in \cite{tb1}--\cite{tb5} that equates Hutchings's contact $1$-form version of embedded contact homology on a given $3$-manifold to a version of this manifold's Seiberg--Witten Floer cohomology. This relationship also has antecedents in the theorem of the second and the third authors proved in \cite{lt} that equates Hutchings's \emph{periodic Floer homology} for fibered $3$-manifolds (see \cite{hs2}) with a version of Seiberg--Witten Floer cohomology.
\par The equivalence between various flavors of Heegaard Floer homology and Seiberg--Witten Floer homology has been conjectured since the discovery of Heegaard Floer homology. See, for example, Conjecture 1.1 in \cite{os2}, Chapter 3.12 in \cite{km} and Conjecture 1 in \cite{kman}.
\begin{ack}
We warmly thank Michael Hutchings for sharing his knowledge of embedded contact homology. We also thank Jonathan Bloom, Tom Mrowka, Peter Ozsv\'{a}th and Zolt\'{a}n Szab\'{o} for helpful thoughts, and we thank Peter Kronheimer and Tom Mrowka for sharing their encyclopedic knowledge of Seiberg--Witten Floer homology before and since their book was published. Last, we thank the staff at the Mathematical Sciences Research Institute (MSRI) for their hospitality. Most of the work reported here was completed during authors' visit to MSRI. The second author was supported in part by the National Science Foundation. The third author was supported in part by a David Eisenbud Fellowship, a Clay Senior Scholar Fellowship and also by the National Science Foundation.
\end{ack}

\section{Heegaard Floer homology of $\m$ and $\ech$ of $\ym$}
\label{s2}
This section describes the $\echinfty$ chain complex on $\ym$ and its relationship with the Heegaard Floer chain complex on $\m$. The first subsection to come summarizes the most salient features of the Heegaard Floer homology related geometry of $\m$. The second subsection describes the $\echinfty$ related geometry on $\ym$. The final subsection explains the precise relationship between the $\echinfty$ chain complex and the relevant Heegaard Floer chain complex. Fix once and for all a $\spinc$ structure on $\m$. This chosen $\spinc$ structure is denoted by $\mathfrak{s}$ in what follows and in the subsequent sections of this article.

\subsection{Heegaard Floer geometry on $\m$}\label{s2.1}
The construction of the three Heegaard Floer homology groups $\hfinfty(\m,\mathfrak{s})$, $\hfminus(\m,\mathfrak{s})$, and $\hfplus(\m,\mathfrak{s})$ requires the data that consists in part of a pair $(f,\psg)$ where $f:\m\rightarrow[0,3]$ is a self-indexing Morse function with a single local maximum and a single local minimum. We take our function $f$ to have $g\geq 1$ index $1$ critical points. It then has $g$ index $2$ critical points and the level set $f^{-1}(3/2)$ is a smooth surface of genus $g$. This surface is denoted by $\Sigma$. What is denoted by $\psg$ is a suitably chosen pseudogradient vector field for $f$ which is smooth on the complement of the critical points of $f$ and is such that $\psg(f)=1$. With the $\spinc$ structure fixed, the pseudogradient vector field $\psg$ is chosen so as to result in what Ozsv\'{a}th and Szab\'{o} call a strongly admissible Heegaard diagram (see Definition 4.10 in \cite{os1}). The diagrams we use are of the kind constructed in the proof of Lemma 5.4 in \cite{os1}.
\par The integral curves of $\psg$ can be used to identify $f^{-1}(1,2)$ with $(1,2)\times\Sigma$ so that the function $f$ appears as the Euclidean coordinate on the $(1,2)$ factor and $\psg$ appears as the corresponding Euclidean vector field. This view of $f^{-1}(1,2)$ led Robert Lipshitz to interpret the Heegaard Floer chain complex as follows (see \cite{rl}). Introduce $\zhf$ to denote the set whose typical element is a collection of $g$ integral curves of $\psg$ that pair the index $1$ and index $2$ critical points of $f$. Any given curve from such a collection starts at an index $1$ critical point of $f$ and ends at an index $2$ critical point of $f$. However, no two curves share the same starting point or the same ending point. The three Heegaard Floer chain complexes are constructed from the free module generated by $\zhf\times\mathbb{Z}$. Any given element in the chain complex can be written as a formal linear combination of elements on $\zhf\times\mathbb{Z}$ with integer coefficients. This free module is referred to in what follows as $\mathbb{Z}(\zhf\times\mathbb{Z})$. 
\par The differential that defines the group $\hfinfty(\m,\mathfrak{s})$ is a certain endomorphism of $\mathbb{Z}(\zhf\times\mathbb{Z})$ whose square is equal to zero. Lipshitz explains in \cite{rl} how to compute the differential from certain geometric data on $\mathbb{R}\times[1,2]\times\Sigma$. The latter requires the choice of a suitably constrained almost complex structure. The almost complex structure is constrained in particular so as to be invariant under the action of $\mathbb{R}$ by translations along the $\mathbb{R}$ factor and to map the Euclidean vector field on the $[1,2]$ factor to that on the $\mathbb{R}$ factor. The almost complex structure is allowed to preserve the tangent plane field to $\Sigma$. However, it need not be constant along the $[1,2]$ factor except in certain prescribed annuli and disks. In particular, it should be constant near the segments that correspond to the integral curves of $\psg$ connecting the index $1$ and index $2$ critical points of $f$. 
\par As explained by Lipshitz in \cite{rl}, the differential for the Heegaard Floer chain complex can be computed using pseudoholomorphic, proper submanifolds with boundary in $\mathbb{R}\times[1,2]\times\Sigma$. These submanifolds exhibit appropriate behavior on the boundary. They also have $2g$ strip-like ends such that the Euclidean coordinate of the $\mathbb{R}$ factor is unbounded from above on $g$ of these ends and it is unbounded from below on the remaining $g$ ends. The asymptotics on these ends with respect to the unbounded Euclidean factor are suitably constrained by data from $\zhf$.
\par Let $\delhf$ denote the differential on $\mathbb{Z}(\zhf\times\mathbb{Z})$. The homology of the resulting chain complex is $\hfinfty(\m,\mathfrak{s})$. As explained by Oszv\'{a}th and Szab\'{o} in \cite{os1}, $\delhf$ preserves the filtration on $\mathbb{Z}(\zhf\times\mathbb{Z})$ defined by the submodules freely generated by $\zhf\times\{i\in\mathbb{Z}\,|\,i\leq k\}$ for any given $k\in\mathbb{Z}$. Granted that such is the case, $\delhf$ defines a differential on the submodule $\mathbb{Z}(\zhf\times\{i\in\mathbb{Z}\,|\,i\leq-1\})$ as well as on the quotient $\mathbb{Z}(\zhf\times\mathbb{Z})/\mathbb{Z}(\zhf\times\{i\in\mathbb{Z}\,|\,i\leq-1\})$. The homology of the former chain complex gives $\hfminus(\m,\mathfrak{s})$, while the homology of the latter chain complex is $\hfplus(\m,\mathfrak{s})$.
\par Any $\spinc$ structure on $\m$ has an associated class in $\mathrm{H}^2(\m;\mathbb{Z})$. The latter is said to be the first Chern class of the $\spinc$ structure. The first Chern class defines via the canonical pairing a homomorphism from the group $\mathrm{H}_2(\m;\mathbb{Z})$ into $\mathbb{Z}$. The image of the latter homomorphism, a subgroup of $\mathbb{Z}$, is generated by an even integer $p\in\mathbb{Z}$. The Heegaard Floer chain complex for the given $\spinc$ structure can be given a relative $\mathbb{Z}/p\mathbb{Z}$ grading so that the differential acts to decrease the grading by $1$.
\par The action of $\mathbb{Z}$ on $\zhf\times\mathbb{Z}$ that translates the $\mathbb{Z}$ factor by $-1$ induces an endomorphism of $\mathbb{Z}(\zhf\times\mathbb{Z})$ that commutes with $\delhf$ and decreases the grading by $2$. The induced endomorphism of $\hfinfty(\m,\mathfrak{s})$, $\hfminus(\m,\mathfrak{s})$, and $\hfplus(\m,\mathfrak{s})$ is the $\umap$-map. As noted in the introduction, there is also a commuting action of $\wedge^\ast(\mathrm{H}_1(\m;\mathbb{Z})/tor)$ on these groups. This action is induced by endomorphisms of the module $\mathbb{Z}(\zhf\times\mathbb{Z})$ that anti-commute with $\delhf$ and decrease the grading by $1$. 

\subsection{Embedded contact geometry on $\ym$}
The manifold $\ym$ is obtained from $\m$ by a surgery that first excises small radius Euclidean balls around the critical points of $f$ and subsequently sutures $g+1$ $1$-handles, copies of $[-1,1]\times\sphere$, to the resulting $2(g+1)$ boundary $2$-spheres. The result is diffeomorphic to $\m\#_{g+1}\handle$. This surgery is done so that the respective spheres around the index $0$ and index $3$ critical points of $f$ connect through one of the attached handles, and the spheres around any given index $1$ critical point connects through a handle to the sphere around one of the index $2$ critical points of $f$. The handle that connects the spheres around the index $0$ and index $3$ critical points of $f$ is denoted by $\ho$. The set of $g$ pairs of index $1$ and index $2$ critical points of $f$ that are defined by the remaining $g$ handles is denoted by $\Lambda$. Given $\pair\in\Lambda$ the corresponding $1$-handle is denoted by $\hp$. The complement in $\y$ of the handles is identified with the complement in $\m$ of the small radius balls around the critical points of $f$. The corresponding parts of $\ym$ and $\m$ are denoted by $\m_{\updelta}$. The manifold $\ym$ is oriented so that the induced orientation on $\m_{\updelta}$ is the opposite of the orientation of $\m$. 
\par The data $(\m,f,\psg)$ are then used to construct a pair $(a,w)$ of a $1$-form and a closed $2$-form on $\ym$ such that $a\wedge w>0$ at all points, and such that $\diff a$ is in the span of $w$. A pair with these properties is said to be a stable Hamiltonian structure. The $2$-form $w$ is chosen so as to be positive on the cross-sectional spheres in $\ho$ and to define an area form on the level sets of $f$ in $\m_{\updelta}$.
\par The $2$-form $w$ is closed and so defines a De Rham cohomology class on $\ym$. To say more about this class, use the Mayer-Vietoris sequence to identify $\mathrm{H}_2(\ym;\mathbb{Z})$ with $$\mathrm{H}_2(\m;\mathbb{Z})\oplus\mathrm{H}_2(\ho;\mathbb{Z})\oplus(\oplus_{\pair\in\Lambda}\mathrm{H}_2(\hp;\mathbb{Z})).$$ 

With this identification understood, the cohomology class of the form $w$ is determined by the linear functional that it defines on the various summands. This linear functional pairs with the generator of $\mathrm{H}_2(\ho;\mathbb{Z})$ that is defined by the cross-sectional spheres in $\ho$ oriented by $w$ to give $2$. It has pairing zero with any $\pair\in\Lambda$ version of the summands $\mathrm{H}_2(\hp;\mathbb{Z})$; and it restricts to the $\mathrm{H}_2(\m;\mathbb{Z})$ summand so as to give the pairing on $\m$ that is defined by the first Chern class of $\mathfrak{s}$. The strong admissibility of the chosen Heegaard diagram on $\m$ is needed to construct a pair $(a,w)$ with all of these listed properties.
\par The kernel of $w$ is generated by a nowhere zero vector field denoted by $\mathit{v}$. The construction of $(a,w)$ is such that $\mathit{v}$ is normal to the cross-sectional spheres in $\ho$ and such that $\mathit{v}=\psg$ on $\m_{\updelta}$. The vector field $\mathit{v}$ is normalized so that $a(\mathit{v})=1$.
\par The embedded contact homology chain complex on $\ym$ is defined using a set $\z$ described as follows. An element $\Theta\in\z$ is a finite set of pairs of the form $(\upgamma,m)$ where $\upgamma$ is a closed integral curve of the vector field $\mathit{v}$ and $m\in\mathbb{Z}$. The collections of pairs in $\z$ with certain constraints on the allowed integer components are used to define embedded contact homology. The precise constraints are described in \cite{hu3}. The relevant subset of $\z$ here is denoted by $\z_{ech}$.
\par The elements in $\z$ are labeled in part by the $\spinc$ structures on $\ym$. The $\spinc$ structure of a given element $\Theta\in\z$ is determined by data consisting of the $2$-plane field $\mathrm{kernel}(a)\subset\mathrm{T}\ym$ and a class that $\Theta$ defines in $\mathrm{H}_1(\ym;\mathbb{Z})$. What follows is the definition of this class. Use $\mathit{v}$ to orient its closed integral curves so as to view them as closed $1$-cycles. If $\upgamma$ is such a curve, use $[\upgamma]$ to denote the corresponding $1$-cycle. The cycle $\sum_{(\upgamma,n)\in\Theta}n[\upgamma]$ is then the desired class in $\mathrm{H}_1(\ym;\mathbb{Z})$. The first Chern class of the corresponding $\spinc$ structure is
\begin{equation}
\label{e2.1}
\mathrm{e}_{\canon}+2\sum_{(\upgamma,n)\in\Theta}n[\upgamma]^{\mathrm{Pd}},
\end{equation}
where $[\upgamma]^{\mathrm{Pd}}$ denotes the class in $\mathrm{H}^2(\ym;\mathbb{Z})$ that is the Poincar\'{e} dual of $[\upgamma]$, and where $\mathrm{e}_{\canon}$ is the Euler class of the oriented $2$-plane field defined by  $\mathrm{kernel}(a)\subset\mathrm{T}\ym$ with the orientation given by $w$. 
\par Our constructions give a natural $1-1$ correspondence between the set of $\spinc$ structures on $\m$ and the set of $\spinc$ structures on $\ym$ that obey the following constraints:
\vspace{11pt}
\begin{itemize}\leftskip -0.35in
\label{e2.2}
\item The first Chern class has pairing $2$ with the cross-sectional spheres in $\ho$.
\item The first Chern class has pairing $0$ with the cross-sectional spheres in $\hp$ for any $\pair\in\Lambda$.
\end{itemize} 
\begin{equation}
\end{equation}
With the preceding understood, of particular interest in what follows is the subset $\z_{ech,\m}\subset\z_{ech}$ of elements whose first Chern class obeys (\ref{e2.2}) and which correspond to the chosen $\spinc$ structure on $\m$. 
\par The pair $(a,w)$ is constructed so as to obtain an essentially explicit description of the set $\z_{ech,\m}$, which we summarize in the next theorem. The latter uses $\textsc{o}$ to denote the set $\{0,1,-1,\{1,-1\}\}$.
\begin{theorem}
\label{t2.1}
The set $\z_{ech,\m}$ is in $1-1$ correspondence with $\zhf\times\prod_{\pair\in\Lambda}(\mathbb{Z}\times\textsc{o})$; and this correspondence is canonical given the choice for $0$ in each $\pair\in\Lambda$ factor of $\mathbb{Z}$. This identification preserves the labeling by $\spinc$ structures.
\end{theorem}
\noindent This theorem is proved in \cite{klt1}. The remarks that follow say something about how the asserted identification comes about. 
\begin{remark}
\label{r1}
The vector field $\mathit{v}$ is normal to the cross-sectional spheres in $\ho$ and as a consequence, the class $\mathrm{e}_{\canon}$ has pairing $2$ with these spheres. It also means that the integral curves of $\mathit{v}$ through $\ho$ have transverse intersections with these spheres with positive local intersection number. These observations with those of the first bullet in (\ref{e2.2}) imply that no curve that contributes to a collection from $\z_{ech,\m}$ can intersect $\ho$.
\end{remark}
\begin{remark}
\label{r2}
The vector field $\mathit{v}$ on $\m_{\updelta}$ is the pseudogradient vector field $\psg$. This implies that any integral curve of $\psg$ that intersects either the $f\leq1$ or the $f\geq2$ part of $\m_{\updelta}$ must cross $\ho$. As a consequence, a given integral curve that appears in an element from $\z_{ech,\m}$ and intersects $\m_{\updelta}$ does so as an integral curve of $\psg$ that is very near an integral curve that connects an index $1$ critical point of $f$ with an index $2$ critical point of $f$. 
\end{remark}
\begin{remark}
\label{r3}
The class $\mathrm{e}_{\canon}$ has pairing $-2$ with the cross-sectional spheres in $\cup_{\pair\in\Lambda}\hp$ with respect to a suitable orientation. This and the second bullet of (\ref{e2.2}) imply that the collection of curves from any given element in $\z_{ech,\m}$ must have intersection number $1$ with the cross-sectional spheres in $\cup_{\pair\in\Lambda}\hp$.
\end{remark}
\begin{remark}
\label{r4}
The geometry in the $1$-handles labeled by $\Lambda$ is such that an integral curve that has intersection number $-1$ with any cross-sectional sphere in these handles will intersect either the $f\leq1$ or the $f\geq2$ part of $\m_{\updelta}$. This together with Remark \ref{r3} has the following consequence. Suppose that $\Theta\in\z_{ech,\m}$. Then, there is precisely one segment in each $\hp$ from the union of the integral curves from $\Theta$ that enters the handle or leaves the handle. Moreover, this segment must cross the handle from the index $2$ critical point end towards the index $1$ critical point end.  
\end{remark}
\begin{remark}
\label{r5}
What is said in Remarks \ref{r2} and \ref{r4} imply that any given element in $\z_{ech,\m}$ defines a canonical corresponding element in $\zhf$. 
\end{remark}
\begin{remark}
\label{r6}
The factor of $\prod_{\pair\in\Lambda}(\mathbb{Z}\times\textsc{o})$ in Theorem \ref{t2.1} labels the possible ways in which a given element in $\zhf$ can be extended over $\cup_{\pair\in\Lambda}\hp$ so as to define an element in $\z_{ech,\m}$. In particular, there are precisely two closed integral curves of $\mathit{v}$ that lie entirely in any given $\hp$. One lies north of the equatorial circle in the central cross-sectional sphere of $\hp$, and the other lies south of this circle. The factor $\textsc{o}$ labels whether a pair from $\Theta$ contains none of these curves, just the one on the northern hemisphere, just the one on the southern hemisphere, or both of these curves. The segment that crosses the handle with positive intersection number does so very near the equator in each cross-sectional sphere. The $\mathbb{Z}$ factor describes the number of times this segment winds around the equator as it traverses the handle. 
\end{remark}
\begin{remark}
\label{r7}
Any given $\Theta\in\z_{ech,\m}$ contains only pairs of the form $(\upgamma,1)$, which follows from the previous remark. This is consistent with Hutchings's constraints on the set of generators of the embedded contact homology chain complex because $\upgamma$'s linearized return map is hyperbolic. 
\end{remark}
\par Any given $\Theta\in\z$ has a \emph{length}, namely, $\sum_{(\upgamma,n)\in\Theta}n\int_{\upgamma}a$. With this in mind, parametrize an element in $\z_{ech,\m}$ by its corresponding $\zhf\times\prod_{\pair\in\Lambda}(\mathbb{Z}\times\textsc{o})$ label. As it turns out, the length of a given element is bounded from below by a fixed multiple of the sum of the absolute values of the integers from the $\prod_{\pair\in\Lambda}(\mathbb{Z}\times\textsc{o})$ factor. This observation motivates the introduction of a filtration of the set $\z_{ech,\m}$ by a nested sequence of finite sets:
\begin{equation}
\label{e2.3}
{\z_{ech,\m}}^1\subset{\z_{ech,\m}}^2\subset\cdots\subset{\z_{ech,\m}}^{\mathrm{L}}\subset\cdots\subset{\z_{ech,\m}}.
\end{equation}
The set ${\z_{ech,\m}}^{\mathrm{L}}$ contains the following sorts of elements. Write a given $\Theta\in\z_{ech,\m}$ using Theorem \ref{t2.1} as $\Theta=(\hat{\upupsilon},(m_{\pair},\mathrm{o}_\pair)_{\pair\in\Lambda})$ with $\hat{\upupsilon}\in\zhf$ and with each $\pair\in\Lambda$ version of $(m_{\pair},\mathrm{o}_\pair)$ denoting a pair in $\mathbb{Z}\times\textsc{o}$. Then, $\Theta\in{\z_{ech,\m}}^{\mathrm{L}}$ if and only if $\sum_{\pair\in\Lambda}(|m_\pair|+2|\mathrm{o}_\pair|)<\mathrm{L}$.

\subsection{The $\echinfty$ chain complex and its homology}\label{ss2.3}
The chain complex for $\echinfty$ is a twisted version of embedded contact homology as described in Section 11 of \cite{hs} and in \cite{hu3}. The definition of twisted embedded contact homology in \cite{hs} requires the choice of a subgroup of $\mathrm{H}_2(\ym;\mathbb{Z})$. The subgroup used here is $\mathrm{H}_2(\m;\mathbb{Z})\oplus(\oplus_{\pair\in\Lambda}\mathrm{H}_2(\hp;\mathbb{Z}))$ seen as a subgroup of $\mathrm{H}_2(\ym;\mathbb{Z})$ via the Mayer--Vietoris sequence. The corresponding chain complex is viewed here as the free module generated by the elements of a certain principal $\mathrm{H}_2(\ho;\mathbb{Z})$-bundle over the discrete set $\z_{ech,\m}$. The latter bundle is denoted in what follows by $\hat{\z}_{ech,\m}$ and the free module that it generates by $\mathbb{Z}(\hat{\z}_{ech,\m})$. Here is the precise definition: a given element in $\mathbb{Z}(\hat{\z}_{ech,\m})$ is a formal, integer weighted sum of finitely many elements from $\hat{\z}_{ech,\m}$. The action of $\mathrm{H}_2(\ho;\mathbb{Z})\simeq\mathbb{Z}$ on $\hat{\z}_{ech,\m}$ induces a corresponding $\mathrm{H}_2(\ho;\mathbb{Z})$ action on the module $\mathbb{Z}(\hat{\z}_{ech,\m})$ and endows it with the structure of a $\mathbb{Z}[\mathpzc{t},\mathpzc{t}^{-1}]$-module. Here, $\mathpzc{t}$ acts as the class represented by the cross-sectional spheres in $\ho$. The module $\mathbb{Z}(\hat{\z}_{ech,\m})$  has a relative embedded contact homology grading by the same cyclic group that grades the corresponding $\zhf$ labels.
\par The differential for the $\echinfty$ chain complex is constructed using a certain endomorphism of $\mathbb{Z}(\hat{\z}_{ech,\m})$ which decreases the relative grading by $1$ and has square equal to zero. The $\umap$-map and the action of $\wedge^{\ast}(\mathrm{H}_1(\m;\mathbb{Z})/tor)$ on the embedded contact homology are likewise given by endomorphisms of $\mathbb{Z}(\hat{\z}_{ech,\m})$. Let $\mathbb{T}$ denote any one of these endomorphisms. The endomorphism $\mathbb{T}$ is defined by its action on generators of $\mathbb{Z}(\hat{\z}_{ech,\m})$. Meanwhile, its action on any given generator can be written as a finite sum
\begin{equation}
\label{e2.4}
\mathbb{T}\hat{\Theta}=\sum_{\hat{\Theta}\in\mathbb{Z}(\hat{\z}_{ech,\m})}\mathpzc{z}_{\hat{\Theta}'\hat{\Theta}}\hat{\Theta}',
\end{equation}
where the coefficients $\mathpzc{z}_{\hat{\Theta}'\hat{\Theta}}$ are integers.
\par The integer coefficients in (\ref{e2.4}) are constructed using pseudoholomorphic submanifolds in $\mathbb{R}\times\ym$ that are defined by a suitably constrained almost complex structure on the latter. The almost complex structure is chosen so as to have various special properties. The most salient features are listed below. The notation uses $s$ to denote the Euclidean coordinate on the $\mathbb{R}$ factor of $\mathbb{R}\times\ym$.
\begin{itemize}\leftskip -0.35in
\label{e2.5}
\item The almost complex structure is invariant under translations along the $\mathbb{R}$ factor of $\mathbb{R}\times\ym$ and it maps $\frac{\partial}{\partial s}$ to $\mathit{v}$.
\item The almost complex structure tames the $2$-form $\diff s\wedge a+w$.
\item The complement in $\ym$ of the union of closed integral curves of $\mathit{v}$ in each $\hp$ is foliated by pseudoholomorphic submanifolds. The leaves intersect $\m_{\updelta}$ as level sets of $f$ and they intersect $\ho$ as the cross-sectional spheres.
\item The almost complex structure is constrained to conform up to the change in orientation to that used by Lipshitz on the relevant part of $\mathbb{R}\times \m_{\updelta}$. 
\end{itemize}
\begin{equation}
\end{equation}
\par The conditions in the first two bullets above are standard requirements for defining embedded contact homology. The conditions in the third bullet is very special to the geometry at hand. In particular, it severely constrains the sorts of pseudoholomorphic submanifolds that contribute to the embedded contact homology differential and other endomorphisms defining additional algebraic structure. The fourth bullet brings the corresponding Heegaard Floer version of $\mathbb{T}$ into the story. These four constraints lead to the characterization of the differential and other endomorphisms given in the upcoming Theorem \ref{t2.3}. Theorem \ref{t2.2} serves to set the stage for Theorem \ref{t2.3}.
\par Theorem \ref{t2.2} makes the formal statement that the relevant endomorphisms can be defined by the rules laid out by Hutchings in \cite{hu3}.
\begin{theorem}
\label{t2.2}
The almost complex structure on $\mathbb{R}\times\ym$ can be chosen to satisfy the conditions in (\ref{e2.5}) to the following end. The differential, the $\umap$-map and the $\wedge^\ast(\mathrm{H}_1(\m;\mathbb{Z})/tor)$ action on $\mathbb{Z}(\hat{\z}_{ech,\m})$ can be defined according to the rules laid out by Hutchings. The latter are represented by endomorphisms that have the form depicted in (\ref{e2.4}). The differential and the action by generators of $\mathrm{H}_1(\m;\mathbb{Z})/tor$ reduce the relative grading by $1$ and anti-commute. The $\umap$-map reduces the relative grading by $2$ and commutes with the other endomorphisms. All of these endomorphisms commute with the action of $\mathrm{H}_2(\ho;\mathbb{Z})$. 
\end{theorem}
\noindent Theorem \ref{t2.2} is proved in \cite{klt1}.
\par The next theorem views the product $\z_{ech,\m}\times\mathbb{Z}$ as a principal $\mathrm{H}_2(\ho;\mathbb{Z})$-bundle over $\z_{ech,\m}$ as follows. The generator given by the cross-sectional sphere in $\mathrm{H}_2(\ho;\mathbb{Z})$ act on the $\mathbb{Z}$ factor so as to send any given integer $k$ to $k-1$. This theorem also uses Theorem \ref{t2.1} to write $\z_{ech,\m}$ as $\zhf\times\prod_{\pair\in\Lambda}(\mathbb{Z}\times\textsc{o})$ and having done that it then moves the $\mathbb{Z}$ factor in $\z_{ech,\m}\times\mathbb{Z}$ to write the latter as $(\zhf\times\mathbb{Z})\times\prod_{\pair\in\Lambda}(\mathbb{Z}\times\textsc{o})$.
\par The last bit of notation concerns conventions. Suppose that $\mathrm{E}$ and $\mathrm{E}'$ are graded chain complexes and that $\Delta$ and $\Delta'$ are respective graded endomorphisms. Then $\mathrm{E}\otimes\mathrm{E}'$ inherits a differential that is written as $\Delta+\Delta'$. These are defined by their actions on the decomposable elements as follows. Suppose that $\mathfrak{e}\in\mathrm{E}$ and $\mathfrak{e}'\in\mathrm{E}'$. Then $(\Delta+\Delta')(\mathfrak{e}\otimes\mathfrak{e}')=\Delta\mathfrak{e}\otimes\mathfrak{e}'+(-1)^{\mathrm{degree}(\Delta')\mathrm{degree}(\mathfrak{e})}\mathfrak{e}\otimes\Delta'\mathfrak{e}'$.
\begin{theorem}
\label{t2.3}
There is a principal $\mathbb{Z}$-bundle isomorphism between $\hat{\z}_{ech,\m}$ and $(\zhf\times\mathbb{Z})\times\prod_{\pair\in\Lambda}(\mathbb{Z}\times\textsc{o})$ with the properties detailed momentarily. Use this isomorphism to identify the $\mathbb{Z}$-module $\mathbb{Z}(\hat{\z}_{ech,\m})$ with $\mathbb{Z}(\zhf\times\mathbb{Z})\otimes\mathbb{Z}(\prod_{\pair\in\Lambda}(\mathbb{Z}\times\textsc{o}))$. This $\mathbb{Z}$-module isomorphism identifies the $\echinfty$ differential as
$$\partial_{ech}=\delhf+\sum_{\pair\in\Lambda}\partial_{\pair},$$
where the endomorphism $\delhf$ denotes the differential on $\mathbb{Z}(\zhf\times\mathbb{Z})$, and each $\partial_\pair$ acts on the corresponding $\mathbb{Z}(\mathbb{Z}\times\textsc{o})$ factor as the square zero endomorphism $\partial_\ast$ given by the rule:
\begin{itemize}\leftskip -0.35in
\item $\partial_{\ast}(m,0)=0$ for each $m\in\mathbb{Z}$,
\item $\partial_{\ast}(m,1)=(m,0)+(m+1,0)$ for each $m\in\mathbb{Z}$,
\item $\partial_{\ast}(m,-1)=(m,0)+(m-1,0)$ for each $m\in\mathbb{Z}$,
\item $\partial_{\ast}(m,\{1,-1\})=(m,-1)-(m,1)+(m+1,-1)-(m-1,1)$ for each $m\in\mathbb{Z}$.
\end{itemize}
Meanwhile, this $\mathbb{Z}$-module isomorphism identifies the $\umap$-map on $\mathbb{Z}(\hat{\z}_{ech,\m})$ with the endomorphism $\umap_{\mathrm{HF}}\otimes\mathbb{I}$ of $\mathbb{Z}(\zhf\times\mathbb{Z})\otimes\mathbb{Z}(\prod_{\pair\in\Lambda}(\mathbb{Z}\times\textsc{o}))$. The isomorphism also identifies the endomorphisms that define the action of $\mathrm{H}_1(\m;\mathbb{Z})/tor$ on $\mathbb{Z}(\hat{\z}_{ech,\m})$ with the endomorphisms that only involve the factor $\mathbb{Z}(\zhf\times\mathbb{Z})$, and they act on this factor so as to define the action of $\mathrm{H}_1(\m;\mathbb{Z})/tor$ on the Heegaard Floer chain complex. 
\end{theorem}
\noindent This theorem is proved in \cite{klt2}.
\par As noted in Section \ref{s2.1}, the differential $\delhf$ on $\mathbb{Z}(\zhf\times\mathbb{Z})$ preserves the submodules $\mathbb{Z}(\zhf\times\{i\in\mathbb{Z}\,|\,i\leq k\})$ defined for each $k\in\mathbb{Z}$. Therefore, Theorem \ref{t2.2} implies that the differential $\partial_{ech}$ on $\mathbb{Z}(\hat{\z}_{ech,\m})$ preserves the filtration of $\mathbb{Z}(\z_{ech}\times\mathbb{Z})$ by the submodules $\mathbb{Z}(\z_{ech}\times\{i\in\mathbb{Z}\,|\,i\leq k\})$ defined for each $k\in\mathbb{Z}$. The incarnation in $\mathbb{Z}(\zhf\times\mathbb{Z})$ of the submodule $\mathbb{Z}(\z_{ech}\times\{i\in\mathbb{Z}\,|\,i\leq 0\})$ is denoted in what follows by $\mathbb{Z}(\hat{\z}^0_{ech,\m})$.
\begin{rmk}
Use Theorem \ref{t2.2} to write $\hat{\z}_{ech,\m}$ as $\z_{ech,\m}\times\mathbb{Z}$. Having done so, Theorem \ref{t2.3} asserts that the $\umap$-map endomorphism acts on a given generator $(\Theta,i)$ to yield $(\Theta,i-1)$. This is not a trivial statement as the embedded contact homology version of the $\umap$-map is defined using pseudoholomorphic submanifolds. By way of comparison, the Heegaard Floer $\umap$-map is defined so as to send any given generator $(\hat{\upupsilon},i)\in\zhf\times\mathbb{Z}$ to $(\hat{\upupsilon},i-1)$.
\end{rmk}
\par The endomorphism $\partial_\ast$ preserves a filtration of $\mathbb{Z}(\mathbb{Z}\times\textsc{o})$ that is defined as follows. Define a function $|\cdot|_{\mathrm{o}}$ on $\textsc{o}$ by the rule $|0|_{\mathrm{o}}=0$, $|-1|_{\mathrm{o}}=1=|1|_{\mathrm{o}}$, and $|\{1,-1\}|_{\mathrm{o}}=2$. For each non-negative integer $\mathrm{L}$, let $\mathcal{V}_{\mathrm{L}}$ denote the submodule of $\mathbb{Z}(\mathbb{Z}\times\textsc{o})$ that is generated by elements of the form $(m,\ast)$ with $|m|+2|\ast|_{\mathrm{o}}<\mathrm{L}$. The nested collection $\{\mathcal{V}_{\mathrm{L}}\}_{\mathrm{L}=0,1,\dots}$ of submodules gives a filtration
\begin{equation}
\label{e2.6}
\mathcal{V}_0\subset\mathcal{V}_1\subset\cdots\subset\mathcal{V}_{\mathrm{L}}\subset\cdots\subset\mathbb{Z}(\mathbb{Z}\times\textsc{o})
\end{equation}
with the property that $\partial_\ast:\mathcal{V}_{\mathrm{L}}\rightarrow\mathcal{V}_{\mathrm{L}}$. The homology of $\partial_{\ast}$ on $\mathbb{Z}(\mathbb{Z}\times\textsc{o})$ is that of the direct limit of the endomorphism that is defined by the restriction of $\partial_{\ast}$ to the nested set of submodules in (\ref{e2.6}). That is to say, any given $\mathrm{z}\in\mathbb{Z}(\mathbb{Z}\times\textsc{o})$ is an element of some $\mathcal{V}_{\mathrm{L}}$, and $\mathrm{z}\in\mathcal{V}_{\mathrm{L}}$ and $\mathrm{z}'\in\mathcal{V}_{\mathrm{L}'}$ which are both in the kernel of $\partial_{\ast}$ represent the same class in the homology of the chain complex $(\mathbb{Z}(\mathbb{Z}\times\textsc{o}),\partial_\ast)$ when there exists $\mathrm{L}''\geq max\{\mathrm{L},\mathrm{L}'\}$ and an element $\mathrm{z}''\in\mathcal{V}_{\mathrm{L}''}$ such that $\mathrm{z}=\mathrm{z}'+\partial_{\ast}\mathrm{z}''$.
\par With these last remarks in mind, reintroduce from (\ref{e2.3}) the filtration $\{{\z_{ech,\m}}^{\mathrm{L}}\}_{\mathrm{L}=1,2,\dots}$ of $\z_{ech,\m}$. Given $\mathrm{L}\in\{1,2,\dots\}$, use ${\hat{\z}_{ech,\m}}^{\hspace{0.3in}\mathrm{L}}\subset\hat{\z}_{ech,\m}$ to denote the corresponding principal bundle over ${\z_{ech,\m}}^{\mathrm{L}}$. The module $\mathbb{Z}({\hat{\z}_{ech,\m}}^{\hspace{0.3in}\mathrm{L}})$ has the filtration
\begin{equation}
\label{e2.7}
\mathbb{Z}({\hat{\z}_{ech,\m}}^{\hspace{0.3in}1})\subset\mathbb{Z}({\hat{\z}_{ech,\m}}^{\hspace{0.3in}2})\subset\cdots\subset\mathbb{Z}({\hat{\z}_{ech,\m}}^{\hspace{0.3in}\mathrm{L}})\subset\cdots\subset\mathbb{Z}(\hat{\z}_{ech,\m}).
\end{equation}
\noindent It follows from Theorem \ref{t2.2} that this filtration is preserved by $\partial_{ech}$, the $\umap$-map, and the action of $\wedge^\ast(\mathrm{H}_1(\m;\mathbb{Z})/tor)$.
\par The next theorem describes the homology of the chain complexes $$(\mathbb{Z}(\hat{\z}_{ech,\m}), \partial_{ech}),\;(\mathpzc{t}\mathbb{Z}(\hat{\z}^0_{ech,\m}), \partial_{ech}),\;\mathrm{and}\;(\mathbb{Z}(\hat{\z}_{ech,\m})/\mathpzc{t}\mathbb{Z}(\hat{\z}^0_{ech,\m}), \partial_{ech}).$$ 
These respective homology groups are denoted by $\echinfty$, $\echminus$, and $\echplus$. Each of the latter has an appropriate relative grading, a corresponding $\umap$-map that reduces the degree by $2$, and an action of $\wedge^\ast(\mathrm{H}_1(\m;\mathbb{Z})/tor)$ whose generators reduce the degree by $1$. In what follows, we use $\vhat$ to denote the graded Abelian group $\mathbb{Z}\oplus\mathbb{Z}$ with the first factor at grading $0$ and the second at grading $1$. The next theorem also refers to a $\umap$-map and an action of $\wedge^\ast(\mathrm{H}_1(\m;\mathbb{Z})/tor)$ on the tensor product of $\hfinfty(\m,\mathfrak{s})$, $\hfminus(\m,\mathfrak{s})$, and $\hfplus(\m,\mathfrak{s})$ with $\vhat^{\otimes g}$. These are defined via the $\umap$-map and the $\wedge^\ast(\mathrm{H}_1(\m;\mathbb{Z})/tor)$ action on $\hfinfty(\m,\mathfrak{s})$, $\hfminus(\m,\mathfrak{s})$, and $\hfplus(\m,\mathfrak{s})$ by simply ignoring the $\vhat^{\otimes g}$ factor.
\begin{theorem}
\label{t2.4}
The $\mathbb{Z}$-module isomorphism depicted in Theorem \ref{t2.3} induces isomorphisms indicated by the vertical arrows in the following commutative diagram
\begin{equation*}\begin{CD}
\llap{\(\cdots  \)}\echminus @>>> \echinfty
@>>>\echplus \rlap{\(\cdots \)}\\
@VVV @VVV @VVV \\
\llap{\(\cdots \)}\hfminus(\m,\mathfrak{s})\otimes\vhat^{\otimes g} @>>>
\hfinfty(\m,\mathfrak{s})\otimes\vhat^{\otimes g}@>>>
\hfplus(\m,\mathfrak{s})\otimes\vhat^{\otimes g}\rlap{\(\cdots  \)}
\end{CD}\end{equation*}
where both the top and the bottom rows are long exact sequences. All homomorphisms preserve the relative gradings, intertwine the respective $\umap$-maps and the $\wedge^\ast(\mathrm{H}_1(\m;\mathbb{Z})/tor)$ actions. Furthermore, the middle vertical arrow intertwines the action of $\mathrm{H}_2(\ho;\mathbb{Z})$.
\end{theorem}
\begin{proof}
The top row is the long exact sequence for the chain complex inclusion $\mathbb{Z}(\hat{\z}^0_{ech,\m})\subset\mathbb{Z}(\hat{\z}_{ech,\m})$ with the differential given by $\partial_{ech}$. The bottom row is the long exact sequence for the short exact sequence defined by the chain complex $\mathbb{Z}(\zhf\times\mathbb{Z})\otimes\vhat^{\otimes g}$ and its subcomplex $\mathbb{Z}(\zhf\times\{i\in\mathbb{Z}\,|\,i\leq k\})\otimes\vhat^{\otimes g}$. The following Lemma is needed to discuss the vertical arrows. The proof of this lemma amounts to a straightforward exercise left to the reader.
\begin{lemma}
\label{l2.5}
The homology of the chain complex $(\mathbb{Z}(\mathbb{Z}\times\textsc{o}),\partial_\ast)$ is isomorphic to $\mathbb{Z}\oplus\mathbb{Z}$. The elements $(0,0)$ and $(0,1)-(1,-1)$ are closed and generate the homology.
\end{lemma}
The isomorphism given by Theorem \ref{t2.3} suggests an iterated double complex spectral sequence to calculate the various embedded contact homology groups. It follows from Lemma \ref{l2.5} that the final non-trivial term is the complex $\mathbb{Z}(\zhf\times\mathbb{Z})\otimes\vhat^{\otimes g}$ with differential $\delhf$. The assertions about the vertical arrows follow directly from this last observation. 
\end{proof}
\section{Seiberg--Witten Floer homology of $\yp$ and $\ech$ of $\ym$}
\label{s3}
This section first describes the relevant versions of the Seiberg--Witten equations on $\yp$ and the corresponding chain complex that computes the associated Seiberg--Witten Floer homology. The last part of the section describes the relationship between the corresponding Seiberg--Witten Floer chain complex and the $\echinfty$ chain complex.
\subsection{The Seiberg--Witten equations on $\yp$}
A detailed discussion of the Seiberg--Witten equations and the corresponding Seiberg--Witten Floer homology groups is given in \cite{km}. What follows is a brief summary of the story for the case at hand. The story here is much like what is told in \cite{tb1}--\cite{tb5} and \cite{lt}.
\par The definition of the Seiberg--Witten equations on any given oriented $3$-manifold requires first the choice of a Riemannian metric. The metric we use to define the Seiberg--Witten equations on $\yp$ is determined using the almost complex structure from (\ref{e2.5}) by the following three rules. First, the vector $\mathit{v}$ has norm $1$. Second, it is orthogonal to the $2$-plane bundle in $\mathrm{T}\yp$ that defines the $+i$ eigenspace in $\mathrm{T}_{\mathbb{C}}\ym$ for the action of the chosen almost complex structure on $\mathrm{T}_{\mathbb{C}}(\mathbb{R}\times\ym)$. Third, the chosen almost complex structure acts as a skew-symmetric endomorphism of this eigenspace. The equations also require the choice of a $\spinc$ structure on $\yp$. 
\par The chosen $\spinc$ structure has an associated rank $2$ complex bundle with Hermitian metric denoted by $\mathbb{S}$. The corresponding Seiberg--Witten equations are equations for a pair $(\mathbb{A},\uppsi)$ where $\mathbb{A}$ is a Hermitian connection on the complex line bundle $det(\mathbb{S}):=\wedge^2\mathbb{S}$ and $\uppsi$ is a smooth section of $\mathbb{S}$. The equations are written in terms of the curvature of $\mathbb{A}$ and the Dirac operator on sections of $\mathbb{S}$ that is defined by the Levi--Civita connection on $\mathrm{T}\yp$ and the connection $\mathbb{A}$ on $det(\mathbb{S})$. To set the notation, introduce $\mathrm{F}_{\mathbb{A}}$ to denote the curvature $2$-form of the connection $\mathbb{A}$ and $\mathcal{D}_{\mathbb{A}}\uppsi$ to denote the action of the Dirac operator on $\uppsi$. Use $\ast$ to denote the Hodge star operator for the chosen Riemannian metric. The equations also refer to a section of $i\mathrm{T}^\ast\yp$ that is defined using $\uppsi$. This $1$-form is written as $\ast\uppsi^{\dagger}\uptau\uppsi$. The definition of $\uppsi^{\dagger}\uptau\uppsi$ is the same that used in, for example, \cite{tb1}. Under certain circumstances, the equations also require a \emph{perturbation term}. The latter is described momentarily. The final input is the choice of a parameter $r\geq1$. Then, the relevant version of the Seiberg--Witten equations read:
\begin{eqnarray}
\label{e3.1}
\nonumber\ast\mathrm{F}_{\mathbb{A}}&=&2r(\uppsi^{\dagger}\uptau\uppsi-i\ast w)+\mathfrak{T}_{(\mathbb{A},\uppsi)}\\\mathcal{D}_{\mathbb{A}}\uppsi&=&\mathfrak{G}_{(\mathbb{A},\uppsi)},
\end{eqnarray}
where $\mathfrak{T}$ and $\mathfrak{G}$ constitute the perturbation term. What is denoted by $\mathfrak{T}_{(\mathbb{A},\uppsi)}$ is a section of $i\mathrm{T}^\ast\yp$ that can be written as $\ast\diff\mathfrak{t}_{(\mathbb{A},\uppsi)}$ where in most cases $\mathfrak{t}_{(\mathbb{A},\uppsi)}$ has non-local dependence on the pair $(\mathbb{A},\uppsi)$. Meanwhile, $\mathfrak{G}_{(\mathbb{A},\uppsi)}$ is a section of $\mathbb{S}$ that also depends in a non-local fashion on $(\mathbb{A},\uppsi)$. These are added to ensure that the set of solutions to (\ref{e3.1}) and also the solutions to the upcoming (\ref{e3.3}) are well-behaved. In any case, these terms have very small norm. The analogous perturbation term is discussed at length in \cite{tb1} in the case when $a$ is a contact form. But for notational changes, the discussion there applies here. This said, the subsequent discussion is worded as if these terms are zero with it understood that everything said applies if they are needed. By way of a parenthetical remark, the term $-2ri\ast w$ on the right hand side of the top equation in (\ref{e3.1}) is an example of what Kronheimer and Mrowka call a \emph{non-exact perturbation}. The Seiberg--Witten equations with non-exact perturbations and their associated Floer homology groups are discussed in Chapters 29-32 of \cite{km} (see Definition 30.1.1). 
\par If $(\mathbb{A},\uppsi)$ is a solution of the equations in (\ref{e3.1}), then so is $(\mathbb{A}-2\gu^{-1}\diff\gu,\gu\uppsi)$ where $\gu$ is any smooth map from $\yp$ into $\mathrm{S}^1$ with $\mathrm{S}^1$ regarded as the unit circle in $\mathbb{C}$. Two solutions that are related in such a way are said to be gauge equivalent. A solution $(\mathbb{A},\uppsi)$ of the equations in (\ref{e3.1}) is said to be \emph{reducible} if $\uppsi$ is identically zero, otherwise it is called \emph{irreducible}.
\par Of particular interest here with regard to our Main Theorem are the versions of (\ref{e3.1}) for the $\spinc$ structure that is used to define $\echinfty$. With this $\spinc$ structure understood, the term $-2ri\ast w$ in (\ref{e3.1}) constitutes what is said to be a \emph{monotone} perturbation. Note that all solutions of the equations in (\ref{e3.1}) for the $\spinc$ structure of interest are irreducible. This is because the first Chern class of $\mathbb{S}$ is represented by the $2$-form $\frac{i}{2\pi}\mathrm{F}_{\mathbb{A}}$ and therefore by $\frac{1}{\pi}rw$ if there were a reducible solution of the equations in (\ref{e3.1}). The latter $2$-form with $r>\pi$ has integral greater than $2$ over any given cross-sectional sphere of $\ho$ and hence the top constraint in (\ref{e2.2}) would be violated.
\subsection{The Seiberg--Witten Floer chain complex on $\yp$}
Let $\go$ denote the subgroup of $\mathrm{C}^{\infty}(\yp;\mathrm{S}^1)$ that is defined as follows. An element of $\go$ represents a class in $\mathrm{H}^1(\yp;\mathbb{Z})$ whose Poincar\'{e} dual in $\mathrm{H}_2(\yp;\mathbb{Z})$ belongs to the group $\mathrm{H}_2(\m;\mathbb{Z})\oplus(\oplus_{\pair\in\Lambda}\mathrm{H}_2(\hp;\mathbb{Z}))$ regarded as a subgroup of the latter via the Mayer--Vietoris sequence. With $r\geq1$ fixed, introduce $\zhsw$ to denote the union over the relevant $\spinc$ structures of the corresponding set of equivalence classes of solutions to the equations in (\ref{e3.1}) with the equivalence relation given by
\begin{equation}
\label{e3.2}
(\mathbb{A},\uppsi)\sim(\mathbb{A}-2\gu^{-1}\diff\gu,\gu\uppsi)
\end{equation}
with $\gu\in\go$. Let $\zsw$ denote the union over the relevant $\spinc$ structures of the corresponding sets of equivalence classes of solutions to the equations in (\ref{e3.1}) with the equivalence relation defined by (\ref{e3.2}) with $\gu\in\mathrm{C}^{\infty}(\y;\mathrm{S}^1)$. The set $\zhsw$ is a principal $\mathrm{H}_2(\ho;\mathrm{Z})$-bundle over $\zsw$. Meanwhile, $\zsw$ is a finite set for a suitably generic choice of a perturbation term to use in (\ref{e3.1}).
\par The Seiberg--Witten Floer chain complex of interest here is the free $\mathbb{Z}$-module generated by the set $\zhsw$. This module is denoted by $\mathbb{Z}(\zhsw)$. The corresponding homology is a version of twisted Seiberg--Witten Floer homology as described in Section 3.7 of \cite{km}. What follows directly lists two properties of the module $\mathbb{Z}(\zhsw)$.
\begin{property}
The module $\mathbb{Z}(\zhsw)$ has a relative grading by the same cyclic group that grades the embedded contact homology chain complex and the Heegaard Floer chain complex.
\end{property}

\begin{property}
The action of $\mathrm{H}_2(\ho;\mathbb{Z})$ on $\zhsw$ gives $\mathbb{Z}(\zhsw)$ the structure of a $\mathbb{Z}[\mathpzc{t},\mathpzc{t}^{-1}]$-module where $\mathpzc{t}$ represents the class defined by the cross-sectional spheres in $\ho$.
\end{property}
The theorem that follows summarizes the salient properties of $\mathbb{Z}(\zhsw)$.
\begin{theorem}
\label{t3.1}
The definition given in \cite{km} for the Seiberg--Witten Floer homology differential supplies a square zero endomorphism, a differential $\delsw$ on $\mathbb{Z}(\zhsw)$. The definitions in \cite{km} also supply an endomorphism that defines the $\umap$-map for $\mathbb{Z}(\zhsw)$ as well as endomorphisms that define the action of $\wedge^{\ast}(\mathrm{H}_1(\m;\mathbb{Z})/tor)$. These endomorphisms have the following properties:
\begin{itemize}\leftskip -0.35in
\item The differential decreases the relative grading by $1$, the $\umap$-map decreases the relative grading by $2$ and the action by generators of $\mathrm{H}_1(\m;\mathbb{Z})/tor$ decrease the grading by $1$. The generators of the $\umap$-map and the action of $\mathrm{H}_1(\m;\mathbb{Z})/tor$ define an action of $\mathbb{Z}[\umap]\otimes\wedge^\ast(\mathrm{H}_1(\m;\mathbb{Z})/tor)$ on the homology of the chain complex.
\item The differential, the $\umap$-map and the action by generators of $\mathrm{H}_1(\m;\mathbb{Z})/tor$ commute with the action of $\mathrm{H}_2(\ho;\mathbb{Z})$.
\item There exists a constant $\upkappa\geq1$ and a section $\mathfrak{Z}\subset\zhsw$ whose significance is described next. If $r\geq\upkappa$, then the $\mathbb{Z}$-module generated by $\zhsw^0:=\cup_{k\geq0}\mathpzc{t}^k\mathfrak{Z}$ is a $\delsw$-invariant submodule of $\mathbb{Z}(\zhsw)$ which is also preserved by the $\umap$-map and the endomorphisms that define the $\wedge^{\ast}(\mathrm{H}_1(\m;\mathbb{Z})/tor)$ action.
\end{itemize}
\end{theorem}
\noindent This theorem is proved in \cite{klt3}. It is assumed implicitly in what follows that the constant $r$ in (\ref{e3.1}) is large enough to invoke all three bullets of Theorem \ref{t3.1}.
\par By way of summary of what is said in \cite{km}, the endomorphisms of $\mathbb{Z}(\zhsw)$ that represent the differential, the $\umap$-map and the action of $\wedge^{\ast}(\mathrm{H}_1(\m;\mathbb{Z})/tor)$ are defined using certain kinds of solutions to the Seiberg--Witten equations on $\mathbb{R}\times\yp$. The latter are equations for an $\mathbb{R}$-dependent pair of a connection on $det(\mathbb{S})$ and a section of $\mathbb{S}$. Let $(\mathbb{A},\uppsi)$ denote such an $\mathbb{R}$-dependent pair. The equations written below use $s$ to denote the Euclidean coordinate on the $\mathbb{R}$ factor of $\mathbb{R}\times\yp$.
\begin{eqnarray}
\label{e3.3}
\nonumber\frac{\partial}{\partial s}\mathbb{A}+\ast\mathrm{F}_{\mathbb{A}}&=&2r(\uppsi^{\dagger}\uptau\uppsi-i\ast w)\\
\frac{\partial}{\partial s}\uppsi+\mathcal{D}_{\mathbb{A}}\uppsi&=&0.
\end{eqnarray}
\par To say something about the kinds of solutions that are relevant, let $(\mathbb{A}_-,\uppsi_-)$ denote a solution to the equations in (\ref{e3.1}). The differential acts on the generator defined by the equivalence class of $(\mathbb{A}_-,\uppsi_-)$ in $\zhsw$ so as to result in a finite linear combination of generators with integer coefficients. Let $(\mathbb{A}_+,\uppsi_+)$ denote a second solution to the equations in (\ref{e3.1}). The latter can represent a generator in this sum only if there exists a solution to the equations in (\ref{e3.3}) with the following properties:
\vspace{11pt}
\begin{itemize}\leftskip -0.35in
\label{e3.4}
\item The map $s\rightarrow(\mathbb{A},\uppsi)|_s$ converges as $s$ tends to $-\infty$ to $(\mathbb{A}_-,\uppsi_-)$.
\item The map $s\rightarrow(\mathbb{A},\uppsi)|_s$ converges as $s$ tends to $\infty$ to $(\mathbb{A}_+-2\gu^{-1}\diff\gu,\gu\uppsi_+)$ with $\gu\in\go$.
\end{itemize}
\begin{equation}
\end{equation}
\par Theorem \ref{t3.1} asserts that the two modules $\mathbb{Z}(\zhsw)$ and $\mathbb{Z}(\zhsw^0)$ together with the differential $\delsw$ define a pair of a chain complex and a subcomplex. Use $\hhbar(\yp)_r$ to denote the homology of the chain complex $(\mathbb{Z}(\zhsw),\delsw)$, $\hhhat(\yp)_r$ to denote the homology of $(\mathbb{Z}(\mathpzc{t}\zhsw^0),\delsw)$, and $\hhcheck(\yp)_r$ to denote the homology of the quotient complex, namely, of the complex $(\mathbb{Z}(\zhsw)/\mathbb{Z}(\mathpzc{t}\zhsw^0),\delsw)$. The next theorem says something about the $r$-dependence of these three homology groups.
\begin{theorem}
\label{t3.2}
Let $\upkappa$ denote the constant from the third bullet of Theorem \ref{t3.1}. Fix $r,r'\geq\upkappa$. There exists a canonical isomorphism between the respective pairs of homology groups \vfill $$(\hhbar(\yp)_r,\hhbar(\yp)_{r'}),\;(\hhhat(\yp)_r,\hhhat(\yp)_{r'}),\;and\;(\hhcheck(\yp)_r,\hhcheck(\yp)_{r'}).$$ These isomorphisms are induced by an $\mathrm{H}_2(\ho;\mathbb{Z})$-equivariant homomorphism from $\mathbb{Z}(\zhsw)$ to $\mathbb{Z}(\hat{\z}_{\mathrm{SW},\y,r'})$ that maps $\mathbb{Z}(\zhsw^0)$ to $\mathbb{Z}(\hat{\z}_{\mathrm{SW},\y,r'}^0)$, preserves their relative gradings, and intertwines the endomorphisms that are used to define the respective differentials, $\umap$-maps and $\wedge^{\ast}(\mathrm{H}_1(\m;\mathbb{Z})/tor)$ actions.
\end{theorem}
\noindent This theorem is also proved in \cite{klt3}.
\par Use Theorem \ref{t3.2} to identify any two $r\geq 1$ versions of the groups $\hhbar(\yp)_r$, $\hhhat(\yp)_r$ and $\hhcheck(\yp)_r$ and so define the $r$-independent groups $\hhbar(\yp)$, $\hhhat(\yp)$ and $\hhcheck(\yp)$.
\subsection{Seiberg--Witten Floer homology on $\yp$ and embedded contact homology on $\ym$}
The upcoming Theorem \ref{t3.3} relates the chain complexes $(\mathbb{Z}(\zhsw),\delsw)$ and $(\mathbb{Z}(\hat{\z}_{ech,\m}),\partial_{ech})$. What follows directly sets up the notation.
\par Theorem \ref{t2.3} describes a certain principal $\mathrm{H}_2(\ho;\mathbb{Z})$-bundle isomorphism from $\hat{\z}_{ech,\m}$ to $\z_{ech,\m}\times\mathbb{Z}$. Subsection \ref{ss2.3} introduces $\mathbb{Z}(\hat{\z}_{ech,\m}^0)\subset\mathbb{Z}(\hat{\z}_{ech,\m})$ to denote the inverse image of $\mathbb{Z}(\z_{ech,\m}\times\{i\in\mathbb{Z}\;|\;i\leq0\})$ via the induced isomorphism of $\mathbb{Z}$-modules. Subsection \ref{ss2.3} also refers to the sets $\{{\z_{ech,\m}}^{\mathrm{L}}\}_{\mathrm{L}=1,2,\dots}$ from (\ref{e2.3}) and the associated filtration of $\mathbb{Z}(\hat{\z}_{ech,\m})$ given by (\ref{e2.7}).
\begin{theorem}
\label{t3.3}
Let $\mathbb{H}_\ast^{\infty}(\yp)$, $\mathbb{H}_\ast^{-}(\yp)$, and $\mathbb{H}_\ast^{+}(\yp)$ denote finitely generated subgroups of the respective groups $\hhbar(\yp)$, $\hhhat(\yp)$, and $\hhcheck(\yp)$. Given these groups, there exists $\mathrm{L}_{\mathbb{H}}\geq 1$, and given $\mathrm{L}\geq\mathrm{L}_{\mathbb{H}}$, there exist $r_{\mathbb{H},\mathrm{L}}\geq\pi$ and $\mathrm{L}'\geq\mathrm{L}$ with the following significance. Take $r\geq r_{\mathbb{H},\mathrm{L}}$ so as to define $\zhsw$. There exists an $\mathrm{H}_2(\ho;\mathbb{Z})$-equivariant, injective map $\hat{\Phi}^r:\hat{\z}_{ech,\m}^{\hspace{0.3in}\mathrm{L}'}\rightarrow\zhsw$ that defines a $\mathbb{Z}$-module monomorphism $$\mathbb{L}^r:\mathbb{Z}(\hat{\z}_{ech,\m}^{\hspace{0.3in}\mathrm{L}'})\rightarrow\mathbb{Z}(\zhsw)$$ with the properties listed below:
\begin{itemize}\leftskip -0.35in
\item $\mathbb{L}^r$ reverses the sign of relative gradings.
\item $\mathbb{L}^r$ induces monomorphisms from $\mathpzc{t}\mathbb{Z}(\hat{\z}_{ech,\m}^{\hspace{0.3in}\mathrm{L}'}\cap\hat{\z}_{ech,\m}^0)$ into $\mathpzc{t}\mathbb{Z}(\zhsw^0)$ and from $\mathbb{Z}(\hat{\z}_{ech,\m}^{\hspace{0.3in}\mathrm{L}'})/\mathpzc{t}\mathbb{Z}(\hat{\z}_{ech,\m}^{\hspace{0.3in}\mathrm{L}'}\cap\hat{\z}_{ech,\m}^0)$ into $\mathbb{Z}(\zhsw)/\mathpzc{t}(\zhsw^0)$.
\item $\mathbb{L}^r$ intertwines $\partial_{ech}$ with $\partial^\ast_{\mathrm{SW}}$, and it also intertwines the endomorphisms that define the respective $\mathbb{Z}[\umap]\otimes\wedge^\ast(\mathrm{H}_1(\y;\mathbb{Z})/tor)$-actions on the $\partial_{ech}$ homology and the $\partial^\ast_{\mathrm{SW}}$ homology.
\item Let ${\mathbb{Q}_{ech}}^{\mathrm{L}}$ denote either one of  $\mathbb{Z}(\hat{\z}_{ech,\m}^{\hspace{0.3in}\mathrm{L}})$, $\mathpzc{t}\mathbb{Z}(\hat{\z}_{ech,\m}^{\hspace{0.3in}\mathrm{L}}\cap\hat{\z}_{ech,\m}^0)$ or\\ $\mathbb{Z}(\hat{\z}_{ech,\m}^{\hspace{0.3in}\mathrm{L}})/\mathpzc{t}\mathbb{Z}(\hat{\z}_{ech,\m}^{\hspace{0.3in}\mathrm{L}}\cap\hat{\z}_{ech,\m}^0)$,  and let ${\mathbb{Q}_{ech}}^{\mathrm{L}'}$ denote the $\mathrm{L}'$ version. Use $\mathbb{Q}_{\mathrm{SW}}$ to denote the corresponding $\mathbb{Z}(\zhsw)$, $\mathpzc{t}\mathbb{Z}(\zhsw^0)$ or \\$\mathbb{Z}(\zhsw)/\mathpzc{t}(\zhsw^0)$ as the case may be. If $\zeta\in{\mathbb{Q}_{ech}}^{\mathrm{L}}$ is such that $\mathbb{L}^{r}(\zeta)=\partial_{\mathrm{SW}}^{\ast}\mathfrak{z}$ for some $\mathfrak{z}\in\mathbb{Q}_{\mathrm{SW}}$, then $\zeta=\partial_{ech}\zeta'$ for some ${\mathbb{Q}_{ech}}^{\mathrm{L}'}$.
\item The subgroups $\mathbb{H}_\ast^{\infty}(\yp)$, $\mathbb{H}_\ast^{-}(\yp)$, and $\mathbb{H}_\ast^{+}(\yp)$ are represented by elements in the respective $\mathbb{L}^r$ images of $\mathbb{Z}(\hat{\z}_{ech,\m}^{\hspace{0.3in}\mathrm{L}})$, $\mathpzc{t}\mathbb{Z}(\hat{\z}_{ech,\m}^{\hspace{0.3in}\mathrm{L}}\cap\hat{\z}_{ech,\m}^0)$ or \\$\mathbb{Z}(\hat{\z}_{ech,\m}^{\hspace{0.3in}\mathrm{L}})/\mathpzc{t}\mathbb{Z}(\hat{\z}_{ech,\m}^{\hspace{0.3in}\mathrm{L}}\cap\hat{\z}_{ech,\m}^0)$.
\end{itemize}
Moreover, if $r,r'\geq r_{\mathbb{H},\mathrm{L}}$, then the homomorphism from Theorem \ref{t3.2} can be chosen to intertwine $\mathbb{L}^r$ and $\mathbb{L}^{r'}$.
\end{theorem}
\noindent This theorem is proved in \cite{klt3}.
\begin{rmk}
The arguments for the proof of Theorem \ref{t3.3} are relatively straightforward modifications of the arguments from \cite{lt} and \cite{tb1}--\cite{tb5}. The most technical and subtle parts of those arguments are in \cite{tb1}--\cite{tb4}; and what is done there is likewise required to prove Theorem \ref{t3.3}. Even so, Theorem \ref{t3.3} needs only the simplest cases from \cite{tb1}--\cite{tb4} because all of the integral curves of $\mathit{v}$ that appear in $\z_{ech,\m}$ do so with integer partner $1$.
\end{rmk}
Theorem \ref{t3.3} leads directly to the following result.
\begin{theorem}
\label{t3.4}
There exist homomorphisms indicated by the arrows in the following commutative diagram
\begin{equation*}
\begin{CD}
\llap{\(\cdots  \)}\echminus @>>> \echinfty
@>>>\echplus \rlap{\(\cdots \)}\\
@VVV @VVV @VVV \\
\llap{\(\cdots \)}\hhhat(\yp) @>>>
\hhbar(\yp) @>>>
\hhcheck(\yp)\rlap{\(\cdots  \)}
\end{CD}
\end{equation*}
such that the vertical arrows are isomorphisms and both rows are long exact sequences. All homomorphisms preserve the relative gradings, intertwine the respective $\umap$-maps and the $\wedge^\ast(\mathrm{H}_1(\m;\mathbb{Z})/tor)$ actions. Furthermore, the middle vertical arrow intertwines the action of $\mathrm{H}_2(\ho;\mathbb{Z})$.
\end{theorem}
\begin{proof}
The horizontal arrows are induced from the respective short exact sequences for the relevant chain complexes. The vertical arrows are induced by the maps $\{\mathbb{L}^r\}_{r\geq1}$ from Theorem \ref{t3.3}.
\par To elaborate, consider the middle vertical arrow. What follows first is its definition. Keep in mind for this purpose that any given class in $\echinfty$ is represented by a closed cycle in some $\mathrm{L}\geq 1$ version of $\hat{\z}_{ech,\m}^{\hspace{0.3in}\mathrm{L}}$. This understood, suppose that $\mathrm{L}\geq1$ and that $\mathrm{z}\in\mathbb{Z}(\hat{\z}_{ech,\m}^{\hspace{0.3in}\mathrm{L}})$ is annihilated by $\partial_{ech}$, and so it represents a class in $\echinfty$. Let $[\mathrm{z}]$ denote this class. Fix $r>>1$ so that $\hat{\Phi}^r$ on $\hat{\z}_{ech,\m}^{\hspace{0.3in}\mathrm{L}'}$ is defined for $\mathrm{L}'\geq\mathrm{L}$. The second and the third bullets of Theorem \ref{t3.3} guarantee that $\mathbb{L}^r\mathrm{z}$ is annihilated by $\delsw$ and hence defines a class in $\hhbar{\yp}$. The second and the third bullets of Theorem \ref{t3.3} imply that $\mathrm{z}$ and $\mathbb{L}^r(\mathrm{z}+\partial_{ech}\mathrm{w})$ define the same class in $\hhbar(\yp)$, and the final assertion of Theorem \ref{t3.3} guarantees that this class is independent of $r$. The homomorphism that is denoted by the middle vertical arrow sends $[\mathrm{z}]\in\echinfty$ to the class represented by $\mathbb{L}^r\mathrm{z}$ in $\hhbar(\yp)$.
\par To see that the middle vertical arrow is injective, let $\mathrm{z}$ be as described above and suppose that $\mathbb{L}^r\mathrm{z}=\delsw \mathfrak{z}$. If $r\geq1$ is sufficiently large and so $\mathrm{L}'$ is sufficiently large, then $\mathfrak{z}$ can be taken to be equal to $\mathbb{L}^r\mathrm{w}$ where $\mathrm{w}\in\hat{\z}_{ech,\m}^{\hspace{0.3in}\mathrm{L}'}$. This is a consequence of the fourth bullet of Theorem \ref{t3.3}. But then $\mathrm{z}=\partial_{ech}\mathrm{w}$ and hence $[\mathrm{z}]=0$.
\par To see that the middle vertical arrow is surjective, suppose that $\mathfrak{z}\in\hhbar(\yp)$ is a given class. If $r\geq1$ is sufficiently large and so $\mathrm{L}'$ is sufficiently large, then the sixth bullet of Theorem \ref{t3.3} finds $\mathrm{z}\in\hat{\z}_{ech,\m}^{\hspace{0.3in}\mathrm{L}'}$ such that $\mathbb{L}^r\mathrm{z}$ represents the class $\mathfrak{z}$. The third and the fourth bullets of Theorem \ref{t3.3} require that $\partial_{ech}\mathrm{z}=0$ and hence $\mathrm{z}$ defines a class in $\echinfty$. As a result, $\mathfrak{z}$ is in the image of the homomorphism indicated by the middle vertical arrow.
\par The arguments just given, after only notational changes, together with the fifth and the sixth bullets of Theorem \ref{t3.3} suffice to prove the assertions for the leftmost and the rightmost vertical arrows.
\end{proof}
\section{Seiberg--Witten Floer homology of $\yp$ and Seiberg--Witten Floer homology of $\m$}
\label{s4}
This section starts by summarizing what is needed about the Seiberg--Witten Floer chain complex on $\m$, and goes on to relate the corresponding three versions of the Seiberg--Witten Floer homology groups to the groups $\hhbar(\yp)$, $\hhhat(\yp)$, and $\hhcheck(\yp)$.
\par To start, fix a smooth Riemannian metric on $\m$ and use the $\spinc$ structure $\mathfrak{s}$ so as to define the Seiberg--Witten equations on $\m$. Use $\mathbb{S}_\m$ to denote the associated rank $2$ Hermitian vector bundle. The relevant versions of the Seiberg--Witten equations involve a pair $(\mathbb{A},\uppsi)$ where $\mathbb{A}$ is a connection on $det(\mathbb{S}_\m)$ and $\uppsi$ is a smooth section of $\mathbb{S}_\m$. The version needed on $\m$ requires the choice of a closed $2$-form that represents the first Chern class of $det(\mathbb{S}_\m)$. Let $\mathrm{c}_{\mathbb{S}_\m}$ denote such a form. The equations read
\begin{eqnarray}
\label{e4.1}
\nonumber\ast\mathrm{F}_{\mathbb{A}}&=&2(\uppsi^{\dagger}\uptau\uppsi-i\pi\ast \mathrm{c}_{\mathbb{S}_\m})+\mathfrak{T}_{(\mathbb{A},\uppsi)}\\
\mathcal{D}_{\mathbb{A}}\uppsi&=&\mathfrak{G}_{(\mathbb{A},\uppsi)},
\end{eqnarray}
As in the case of (\ref{e3.1}), the perturbation terms are only needed to ensure that the solution set of these equations and the corresponding $\mathbb{R}\times\m$ version of (\ref{e3.3}) are well behaved. As was done in Section \ref{s3}, the subsequent discussion is phrased as if they are absent.
\par The equations lead to what Kronheimer and Mrowka call \emph{balanced} Floer homology groups in \cite{km}. These are discussed in Chapters 29 and 30 of \cite{km} (see Definition 30.1.1). There are three versions of balanced Seiberg--Witten Floer homology groups with coefficient ring $\mathbb{Z}$. These are the groups $\hmbar(\m,\mathfrak{s},c_b)$, $\hmhat(\m,\mathfrak{s},c_b)$, and $\hmcheck(\m,\mathfrak{s},c_b)$ that appear in our Main Theorem. They reduce to the standard versions of Seiberg--Witten Floer homology when the first Chern class of $det(\mathbb{S}_\m)$ is torsion, otherwise they are different. For example, $\hmbar(\m,\mathfrak{s}_\m)$ is zero when the associated first Chern class is non-torsion while $\hmbar(\m,\mathfrak{s}_\m,c_b)$ is not zero in general. Each of these groups has a relative grading, a $\umap$-map that decreases the relative grading by $2$, and a commuting action of $\wedge^{\ast}(\mathrm{H}_1(\m;\mathbb{Z})/tor)$. These three groups fit into a long exact sequence whose homomorphisms intertwine the $\umap$-maps and the $\wedge^{\ast}(\mathrm{H}_1(\m;\mathbb{Z})/tor)$ actions.
\par The manifold $\yp$ is obtained from $\m$ by connected summing the latter with $g+1$ copies of $\handle$. One of these copies accounts for the handle $\ho$ and the others account for the various $\pair\in\Lambda$ versions of $\hp$. This said, the three balanced versions of the Seiberg--Witten Floer homology on $\m$ are related to the corresponding three versions that appear in Theorem \ref{t3.4} using Seiberg--Witten Floer homology connected sum theorems. A $g$-fold iteration of one such theorem deals with the $\pair\in\Lambda$ labeled copies of $\handle$. A second sort of connected sum theorem deals with the $\ho$ labeled copy. The connected sum theorems lead directly to
\begin{theorem}
\label{t4.1}
There exist homomorphisms indicated by the arrows in the following commutative diagram
\begin{equation*}\begin{CD}
\llap{\(\cdots  \)}\hhhat(\yp) @>>>\hhbar(\yp)
@>>>\hhcheck(\yp) \rlap{\(\cdots \)}\\
@VVV @VVV @VVV \\
\llap{\(\cdots \)}\hmhat(\m,\mathfrak{s},c_b)\otimes\vhat^{\otimes g} @>>>
\hmbar(\m,\mathfrak{s},c_b)\otimes\vhat^{\otimes g}@>>>
\hmcheck(\m,\mathfrak{s},c_b)\otimes\vhat^{\otimes g}\rlap{\(\cdots  \)}
\end{CD}\end{equation*}
such that the vertical arrows are isomorphisms and both rows are long exact sequences. All homomorphisms preserve the relative gradings, intertwine the respective $\umap$-maps and the $\wedge^\ast(\mathrm{H}_1(\m;\mathbb{Z})/tor)$ actions.
\end{theorem}
\noindent A proof of this theorem will appear in \cite{klt4}.
\par The Main Theorem from the introduction is an immediate consequence of Theorems \ref{t2.4}, \ref{t3.4} and \ref{t4.1}.
\begin{rmk}
The use of two different connected sum theorems owes allegiance to the distinction in (\ref{e2.2}). The connected sum theorem that deals with the $\pair\in\Lambda$ labeled copies of $\handle$ is the simpler of the two. This theorem accounts for the $g$ factors of $\vhat$ that appear in the bottom row of the diagram in Theorem \ref{t4.1}. Suffice it to say for now that any given factor of $\vhat$ reflects the following geometric fact. The space of gauge equivalence classes of flat $\mathrm{U}(1)$ connections on $\handle$ is a manifold, $\mathrm{S}^1$, whose cohomology is isomorphic to $\vhat$. The second sort of connected sum theorem is a very much more subtle affair. Its proof in a slightly different context is outlined in \cite{yjl}. As explained in \cite{yjl}, the theorem owes allegiance in part to various results by J. D. S. Jones regarding $\mathrm{S}^1$-equivariant homology (see \cite{j}). As noted at the end of the Introduction, the existence of this second sort of connected sum theorem was predicted by Mrowka and Oszv\'{a}th; and Bloom, Mrowka and Oszv\'{a}th have a proof of a closely related version \cite{bmo}.
\end{rmk}

\end{document}